\documentclass[12pt,reqno]{amsart}
\usepackage{amssymb,latexsym}

\newcommand{\Z}{{\mathbb Z}}

\theoremstyle{plain}
\newtheorem{thm}{Theorem}

\begin{document}

\title{Compositions, Partitions, and Fibonacci Numbers}
\author{Andrew V. Sills}
\address{Department of Mathematical Sciences\\
                Georgia Southern  University\\
                Statesboro, Georgia\\
                30460-8093 USA}
\email{asills@georgiasouthern.edu}
\begin{abstract}
A bijective proof is given for the following theorem:
the number of compositions of $n$ into odd parts equals the number of compositions of $n+1$
into parts greater than one.   Some commentary about the history of partitions and compositions
is provided.
\end{abstract}

\maketitle

\section{Introduction} 

   A \emph{composition} of an integer $n$ is a representation of $n$ as a sum of
 positive integers, for example the eight compositions of $4$ are as follows:
 \[ 4, \quad  3+1,\quad 1+3,\quad 2+2, \quad 2+1+1, \quad1+2+1, \quad1+1+2,\quad 1+1+1+1. \]
  A \emph{partition} of $n$ is a representation of $n$ as a sum of positive integers where
  the order of
 the summands is considered irrelevant.  Thus $2+1+1$, $1+2+1$ and $1+1+2$ are
 three distinct compositions, but are all considered to be the same partition of $5$.
 
  The individual summands of a composition or partition are called its \emph{parts}.
  By convention, the parts of a partition are written in weakly decreasing order.
  Thus the five partitions of $4$ are
  \[ 4, \quad  3+1,\quad 2+2, \quad 2+1+1, \quad 1+1+1+1. \]
 
 A composition of $n$ with $\ell$ parts may be represented graphically by $n$ unit lengths
separated by $\ell-1$ nodes in such a way as to depict a part $j$ by $j$ adjacent unit
lengths bounded by a node on either side (in the case of an interior part) and
a node on one side (in the case of the first or last part).
See MacMahon~\cite[Sec. IV, Ch. 1, p. 151, \P 125] {MacMahon}).
Let us agree to call this representation the \emph{MacMahon graph} of a composition.

For example, the composition $2+4+1+1+5$ of $13$ has the following MacMahon graph:
\begin{equation*} \label{graph24115} - - \cdot - - - - \cdot - \cdot - \cdot - - - - - \end{equation*}
It is then straightforward to encode the MacMahon graph of a composition of $n$ with $\ell$ parts
as a bit sequence of length $n-1$
consisting of $\ell-1$ ones and $n-\ell$ zeros, where as we read
the graph from left to right, we choose a $1$ if a node is present between two unit lenghts,
and a $0$ if no node is present.  For consistency, we may as well call this the
\emph{MacMahon bit sequence} of a composition.
  
  Thus the MacMahon bit sequence corresponding to the composition $2+4+1+1+5$ is 
\[ \texttt{01 0001 1 1 0000}. \]
 We have chosen to put a space after each ``1" in the bit sequence, as this makes it easy to read off the corresponding
 composition.

 It follows from the MacMahon bit sequence that there are $2^{n-1}$ compositions of $n$,
 since the bit sequence is of length $n-1$, and each bit may take on either of two values. 
 Thus, for $n$ a positive integer, we have the formula
 \begin{equation}\label{cofn}
 c(n) = 2^{n-1}
 \end{equation}
 where $c(n)$ denotes the number of compositions of $n$.
 
 The theory of partitions began with Euler in the mid-eighteenth century.
 In an effort to understand certain aspects of partitions, Euler introduced the idea of a \emph{generating function} of a sequence $\{ a_n \}$, that is,
 he encoded the sequence as the coefficients of a power series $\sum_{n=0}^\infty a_n x^n.$
 In particular, Euler showed that the generating function
 for $p(n)$, the number of partitions of $n$, can be expressed as an elegant infinite
 product:
 \begin{equation}\sum_{n=0}^\infty p(n) x^n = \prod_{j=1}^\infty \frac{1}{1-x^j},
 \label{PartitionGF}\end{equation}
 where $|x|<1$ to ensure convergence if $x$ is taken to be a complex variable. 
 This is necessary when analytic properties of~\eqref{PartitionGF} are studied. 
 When~\eqref{PartitionGF} is used for combinatorial purposes, $x$ may be taken to be
 a formal variable.  Euler also showed that~\eqref{PartitionGF} implies the following
 recurrence for $p(n)$:
 \begin{equation} \label{EulerRec} p(n) = \sum_{j=1}^\infty (-1)^{j-1} \left( p\left(n- \frac{j(3j-1)}{2} \right) + p\left(n- \frac{j(3j+1)}{2} \right)\right)
 . \end{equation}
  
 We pause to remark that in contrast to the extremely simple formula~\eqref{cofn} for $c(n)$, the
behavior of $p(n)$ is much more complicated.  In 1937, Hans Rademacher~\cite{R37}, building
on earlier work by Hardy and Ramanujan~\cite{HR18}, proved the following formula for $p(n)$:
\begin{equation}\label{Rpofn}
 p(n) = \frac{1}{\pi\sqrt{2}}\sum_{k=1}^\infty\sqrt{k} A_k(n) \frac{d}{dn}  \left\{ 
 \frac{\sinh\left( \frac{\pi}{k} \sqrt{\frac 23 \left( n-\frac{1}{24}\right)} \right)}{\sqrt{n-\frac{1}{24}}} 
 \right \},
\end{equation}
where $A_k(n)$ is a Kloosterman-type sum, 
\[ A_k(n) = \underset{\gcd(h,k)=1}{\sum_{0\leq h < k}} \exp\left( \pi i s(h,k) - \frac{2\pi i n h}{k} \right) ,\]
and $s(h,k)$ is a Dedekind sum,
\[  s(h,k) = \sum_{j=1}^{k-1} \left(\left(  \frac j k \right)\right) \left( \left( \frac{hj}{k}  \right) \right), \]
with \[ ((x)):= 
\left\{ \begin{array}{ll} 
x - \lfloor x \rfloor - \frac{1}{2}, &\mbox{if $x\not\in\Z$} \\
0, & \mbox{if $x\in \Z$} \end{array} \right.
 .\]
 Peter Shiu adapted~\eqref{Rpofn} to compute $p(10^6)$.  (See~\cite{PS})
 
 Eq.~\eqref{Rpofn} was the only exact explicit formula known for $p(n)$ until Ken Ono and Jan
 Bruinier amazed the mathematical community in early 2011 with the announcement~\cite{BO}
 of a new formula for $p(n)$ as a \emph{finite} sum of algebraic numbers.
 The Ono-Bruinier formula for $p(n)$, however, is by no means elementary: the algebraic numbers
 in question are singular moduli for a $\Gamma_0(6)$ weight $-2$ meromorphic modular
 form expressible in terms of the quasimodular Eisenstein series 
 \[ E_2(z):= 1 - 24 \sum_{n=1}^\infty \sum_{d \mid n} d q^n \] and
 the Dedekind eta function
 \[ \eta(z) := e^{\pi i z/12} \prod_{n=1}^\infty (1 - e^{2\pi i z n}). \] 

Often one wishes to consider sets of partitions where there is
some restriction on which parts may appear.
Perhaps the most famous identity in the theory of partitions is
\begin{thm}[Euler's partition theorem] \label{EulerOddDist} The number of partitions of $n$ into odd parts equals the number of partitions of
$n$ into distinct parts.
\end{thm}

For example, consider the partitions of $8$ into distinct parts:
\[ 8 \qquad 7+1 \qquad 6+2 \qquad 5+3 \qquad 5+2+1 \qquad 4+3+1. \]
There are six in all.
Now consider the partitions of $8$ into odd parts:
\[ 7+1 \quad 5+3 \quad 5+1+1+1 \quad 3+3+1+1 \quad 3+1+1+1+1+1 \quad 1+1+1+1+1+1+1+1 \]
Again, there are six in all.

Notice that while Theorem~\ref{EulerOddDist} predicts that there are the same number of
partitions of $8$ into distinct parts as there are partitions of $8$ into odd parts, the fact that
there happen to be six such partitions in each class is \emph{not} predicted by Theorem~\ref{EulerOddDist}.

Let $q(n)$ denote the number of partitions of $n$ into odd parts, or into distinct parts.  THen
$q(n)$ may be computed via the recurrence
\[ q(n) + \sum_{k=1}^\infty (-1)^k \left( q\left( n - 3k+1) \right) + q\left( n - 3k-1 \right)   \right)
 = \left\{ \begin{array}{ll} 1 & \mbox{if $n=m(m+1)/2$} \\ 0 &\mbox{otherwise}  \end{array} \right.,\]
which is Theorem 1 in~\cite{ORW}.

Euler proved Theorem~\ref{EulerOddDist} in a book published in 1748 \cite[p. 275]{Euler} using generating functions.  The first bijective proof 
of Theorem~\ref{EulerOddDist} was given by J. W. L. Glaisher~\cite{Glaisher} in 1883.

  For a thorough yet gentle exposition of these two proofs of Theorem~\ref{EulerOddDist}, please see Chapters 5
and 2 respectively of Andrews and Eriksson's book \emph{Integer Partitions}~\cite{AE04}.

It would seem natural to seek an analogous identity involving compositions of
$n$ with odd parts.

We offer:
\begin{thm} \label{MainThm} The number of compositions of $n$ into odd parts equals the number of compositions 
of $n+1$ into parts greater than one.
\end{thm}

Theorem~\ref{MainThm} is an immediate consequence of two older results:
Cayley~\cite{Cayley} showed that the number of compositions of $n+1$ into parts greater than one equals 
the $n$th Fibonacci number $F_{n}.$
That the number of compositions of $n$ into odd parts 
equals $F_n$ is observed a number of places in the literature.  The earliest reference this author found was
in Volume 1 of Richard Stanley's \emph{Enumerative Combinatorics}~\cite[p. 46, ex. 14]{Stanley}, which
was published in 1986.   However, it seems likely that the result is much older.

Here we define
$F_n$  by
 \[ F_0=0,\qquad F_1=1\] and
 \[ F_{n}=F_{n-1}+F_{n-2} \qquad\mbox{when $n\geq 2$}. \]

Our objective in the next section is to provide a bijective proof of Theorem~\ref{MainThm}.
But before doing so, we need to introduce the conjugate of a composition.

 The \emph{conjugate} $\mathbf{c'}$ of a composition $\mathbf{c}$ is the composition obtained from $\mathbf{c}$
by taking the bit complement of the MacMahon bit sequence of $\mathbf{c}$, i.e. change all
of the zeros to ones and vice versa.  
  Thus the conjugate of the composition $2+4+1+1+5$ is 
  $1+2+1+1+4+1+1+1+1$ 
because the MacMahon bit sequence of the former is 
 \[ \texttt{ 0100 0111  0000} \]
while that of the latter is
 \[ \texttt{ 1011 1000 1111}. \]
 
 \section{Bijective proof of Theorem~\ref{MainThm} with a guiding example}
 
 Let us begin with a composition $\mathbf{a}$ of $n$ into $\ell$ parts,
 $a_1 + a_2 + \cdots + a_\ell$ in which each each part $a_i$ is odd.
We wish to map $\mathbf{a}$ to a composition of $n+1$ in which all parts are
greater than $1$.

As we proceed, let us visualize an example.  Let us take 
$\mathbf{a}$ to be $1+1+1+9+1+1+5+3$,
which has MacMahon bit sequence
\verb+111 000000001 1 1 00001 00+.

Notice that because all of the parts in $\mathbf{a}$ are odd, the corresponding
MacMahon bit sequence must have zeros appear in strings of even length.

Let us now map $\mathbf{a}$ to its conjugate composition $\mathbf{a'}$, which
is a composition of $n$ of length $n- \ell + 1$.  The number of parts in
$\mathbf{a'}$ must be odd as $n$ and $\ell$ must be of the same parity
from the elementary fact that the sum of an even number (resp. odd number)
of odd integers is even (resp. odd).  

Recalling that the MacMahon bit sequence of $\mathbf{a}$ has all of
its zeros appearing in strings of even length,
the MacMahon bit sequence for $\mathbf{a'}$ must have the property
that all of its ones appear in strings of even length.
In our example, we have the MacMahon bit sequence of $\mathbf{a'}$ as
\verb+0001 1 1 1 1 1 1 1 0001 1  1 1 01 1+, so in our example $\mathbf{a'}$ is
\[ 4+1+1+1+1+1+1+1+ 4+1+ 1+1+ 2+ 1+1. \]

That the bit sequence of $\mathbf{a'}$ consists of pairs of ones (sometimes
separated by some zeros) means that every even index part, i.e. the
second, fourth, sixth, etc. part, must be a $1$.

Thus let us map $\mathbf{a'}$ to the composition $\mathbf{b}$ obtained
from $a'_1 + a'_2 + \cdots + a'_{n-\ell+1}$ by summing adjacent pairs of parts, i.e.
let \[ b_i = a'_{2i-1} + a'_{2i} \]
for $i=1,2, \dots, \frac{n-\ell}{2};$
and let the last part of $\mathbf{b}$ equal the last part of $\mathbf{a'}$:
\[ b_{\frac{n-\ell}{2} + 1 } = a'_{n-\ell+1}. \]

So in our example, we have $\mathbf{b}$ given by
\[ 5+2+2+2+5+2+3+1. \]

Thus we see that $\mathbf{b}$ is \emph{almost} our target; it is a composition of $n$
where all parts other than the last part are greater than $1$.  But this is easily fixed.
Simply map $\mathbf{b}$ to the composition $\mathbf{c}$ obtained from $\mathbf{b}$
by increasing the last part of $\mathbf{b}$ by one.   Thus, $\mathbf{c}$ is a 
composition of $n+1$ in which all parts are greater than one, as desired. 
In our example, we have that the image of $1+1+1+9+1+1+5+3$ is
$5+2+2+2+5+2+3+2$.

The steps are easily reversible: begin with a composition of $n+1$ into parts
greater than one, reduce the last part by $1$, split each part $j$ (other than the
last part) into the pair of parts $j-1, 1$, and conjugate the resulting 
composition.  The final composition
is a composition of $n$ with all parts odd.
\qed

\section{Commentary on enumeration}
Theorem~\ref{MainThm} tells us that the number of compositions of $n$ into odd parts equals the
number of compositions of $n+1$ into parts greater than one, but does not indicate how many
such compositions there are for a given $n$.  Likewise, Euler's Theorem~\ref{EulerOddDist}
indicates that there are as many partitions of $n$ into odd parts as there are partitions of $n$
into distinct parts, but does not indicate how many such partitions there are for a given $n$.

  The first to find an exact explicit formula for the number $q(n)$ of partitions of $n$ into
odd parts was Loo-Keng Hua~\cite[p. 195]{Hua}.  Later, Peter Hagis~\cite{Hagis} gave the following
Rademacher-type convergent series representation for $q(n)$. 
\begin{equation}
  q(n) = \frac{\pi}{\sqrt{24n+1}} \underset{2\nmid k}{\sum_{k\geq 1}} 
  k^{-1} \underset{\gcd(h,k)=1}{\sum_{1\leq h < k}}
   \exp\Bigg( {\pi i \left( t(h,k)- \frac{2 n h}{k} \right)} \Bigg) I_1 \left( \frac{\pi\sqrt{48n+2}}{12k} \right),
\end{equation}
where \[ t(h,k) = \sum_{j=1}^k \left(\left( \frac{2j-1}{2k} \right)\right) \left(\left( \frac{h(2j-1)}{k} \right) \right), \]
and once again \[ ((x)):= 
\left\{ \begin{array}{ll} 
x - \lfloor x \rfloor - \frac{1}{2}, &\mbox{if $x\not\in\Z$} \\
0, & \mbox{if $x\in \Z$} \end{array} \right.
 ,\] and $I_1(z)$ is the modified
Bessel function of order 1~\cite[p. 77, Eq. (2)]{Watson}.

  In contrast, letting $Q(n)$ denote the number of compositions into odd parts, we
have the very simple formula
\begin{equation} \label{QFib} Q(n) = F_n. 
\end{equation}

It is an easy exercise to establish Eq.~\eqref{QFib} combinatorially, as in~\cite{Grimaldi}.
Of course, Binet's formula immediately implies 
\[ Q(n) =\frac{\left( 1+\sqrt{5} \right)^n - \left( 1-\sqrt{5}\right)^n}{2^n \sqrt{5}}, \]
giving us a direct, non-recursive formula for $Q(n)$ that is much simpler than anything known for 
$q(n)$.  However, this formula for $Q(n)$ is not useful for large values of $n$ because of the
roundoff error generated in computing $(1\pm \sqrt{5})^n$.

  The theory of partitions has a longer and more varied literature than that of compositions.
Perhaps this is due to mathematicians feeling that the theory of partitions is deeper, and
hence more interesting than that of compositions.  After all, the generating function for
$p(n)$ is a certain infinite product which (up to a trivial multiple) just so happens to be
a modular form.  
On the other hand, the generating function
for $c(n)$ can easily be seen, using ~\eqref{cofn} and then summing the geometric series
to be 
\[ \sum_{n=1}^\infty c(n) x^n = \sum_{n=1}^\infty 2^{n-1} x^n = \frac{x}{1-2x}, \]
a mere rational function.

  But before we too hastily dismiss compositions as the less worthy 
 relatives of partitions,
let us turn back to Euler's partition theorem (Theorem~\ref{EulerOddDist}) for some further inspiration.   Because Theorem~\ref{EulerOddDist}  involves partitions into
odd parts, we considered compositions into odd parts.  But Theorem~\ref{EulerOddDist}
also involves partitions into distinct parts.  So, what about \emph{compositions} into
distinct parts?  These are not so easily dealt with.  

  We can make a start by observing that for any partition into $\ell$ distinct parts, there
corresponds $\ell !$ compositions into $\ell$ distinct parts.  Using ideas that follow from 
Euler's, it can be shown that the generating function for partitions into $\ell$ distinct
parts is
  \[ \frac{x^{\ell(\ell+1)/2} }{(1-x)(1-x^2)(1-x^3)\cdots(1-x^\ell) }, \]
and thus the generating function for compositions into $\ell$ distinct parts is
   \[ \frac{x^{\ell(\ell+1)/2}\ \ell!  }{(1-x)(1-x^2)(1-x^3)\cdots(1-x^\ell) }, \]
and thus the generating function for compositions into distinct parts is
  \[ \sum_{\ell=0}^\infty \frac{x^{\ell(\ell+1)/2}\ \ell!  }{(1-x)(1-x^2)(1-x^3)\cdots(1-x^\ell) }, \]
a far cry from the simplicity of the generating function for $c(n)$.
In fact, in a 1995 paper~\cite{RK} entitled ``Compositions with distinct parts," 
Bruce Richmond and Arnold Knopfmacher remark that their ``analysis is more 
complicated than is usual for compositions problems.  The results imply however
that the behaviour of these functions is of comparable complexity to partition problems."

\section*{Acknowledgements}
The author thanks George Andrews for suggesting that he look for a bijective proof of 
Theorem~\ref{MainThm}.  The author also thanks the referee for making a number of
helpful suggestions.

\medskip

\noindent AMS Classification Numbers: 11B39, 05A17


\begin{thebibliography}{99}

\bibitem{AE04} G. E. Andrews and K. Eriksson, \emph{Integer Partitions}, Cambridge University Press, 2004.
\bibitem{BO} J. Bruiner and K. Ono, Algebraic formulas for the coefficients of
half-integral weight harmonic weak Maass forms, preprint, 2011.
\texttt{http://www.aimath.org/news/partition/brunier-ono}
\bibitem{Cayley} A. Cayley, Theorems in trigonometry and on partitions,
 \emph{Collected Mathematical Papers}, vol. 10, 16.
\bibitem{Euler} L. Euler, \emph{Introduction to Analysis of the Infinite: Book I}, 
translation of \emph{Introductio in Analysin Infinitorum} (1748) to English from the original Latin by J. D. Blanton, Springer-Verlag,
1988.
\bibitem{Glaisher} J. W. L. Glaisher, A theorem in partitions, \emph{Messenger of Math.}
\textbf{12} (1883), 158--170.
\bibitem{Grimaldi} R. P. Grimaldi, Compositions with odd summands, 
\emph{Congressus numerantium} \textbf{142} (2000), 113--127. 
\bibitem{Hagis} P. Hagis, Partitions into odd summands, 
\emph{Amer. J. Math.} \textbf{85} (1963),
213--222.

\bibitem{HR18} G. H. Hardy and S. Ramanujan, 
Asymptotic formulae in combinatory analysis,
\emph{Proc. London Math Soc.} (2) \textbf{17} (1918) 75--115.

\bibitem{Hua} L.-K. Hua, On the partitions of a number into unequal parts,
\emph{Trans. Amer. Math. Soc.} \textbf{51} (1942), 194--201.

\bibitem{MacMahon} P. A. MacMahon, \emph{Combinatory Analysis}, vol. 1, Cambridge: at the University
Press, 1915.

\bibitem{R37} H. Rademacher, On the partition function $p(n)$, \emph{Proc. London Math
Soc.} (2) \textbf{37} (1937) 241--254.

\bibitem{RK} B. Richmond and A. Knopfmacher, Compositions with distinct parts,
 {\it Aequationes Math. {\bf 49}} (1995), 86-97.
 
 \bibitem{ORW} K. Ono, N. Robbins, and B. Wilson, Some recurrences for arithmetical functions,
 \emph{J. Indian Math. Soc.} \textbf{62} (1996) 29--51.
 
 \bibitem{PS} P. Shiu, {Computations of the partition function}  \emph{Math. Gazette} \textbf{81}
 (1997) 45--52.

\bibitem{Stanley} R. P. Stanley, \emph{Enumerative Combinatorics}, vol. 1, Wadsworth, 1986; reissued,
Cambridge, 1999.

\bibitem{Watson} G. N. Watson, \emph{A Treatise on the Theory of Bessel Functions}, 2nd ed.,
Cambridge Univ. Press, 1966.


\end{thebibliography}
\end{document}